# INFINITESIMAL CONTRACTION OF THE FEIGENBAUM RENORMALIZATION OPERATOR IN THE HORIZONTAL DIRECTION

DANIEL SMANIA

ABSTRACT. We describe a new proof of the exponential contraction of the Feigenbaum renormalization operator in the hybrid class of the Feigenbaum fixed point. The proof uses the non existence of invariant line fields in the Feigenbaum tower (C. McMullen), the topological convergence (D. Sullivan), and a new infinitesimal argument, different from previous methods by C. McMullen and M. Lyubich.

1. INTRODUCTION AND STATEMENT OF THE RESULT.

Let $f\colon U_0 \to V_0$ be a quadratic-like map. This means that $f$ is a ramified holomorphic covering map of degree two, where $U$ and $V$ are simply connected domains, $U \Subset V$. We also assume that the filled-in Julia set of $f$, $K(f) := \cap_n f^{-n} V$, is connected. We say that $f$ is renormalizable with period two if there exist simply connected subdomains $U_1$, $V_1$ so that $f^2\colon U_1 \to V_1$ is also a quadratic-like map with connected filled-in Julia set.

Two quadratic like maps $f\colon U_f \to V_f$ and $g\colon U_g \to V_g$, both with connected filled-in Julia set, defines the same quadratic-like germ if $K(f)$ coincides with $K(g)$ and $f$ coincides with $g$ in a neighborhood of $K(f)$. If $f$ is renormalizable, then the renormalization of the germ defined by $f$ is the unique quadratic-like germ defined by the normalization of any possible induced map $f^2\colon U_1 \to V_1$ which are quadratic-like maps with connected filled-in Julia set (normalize the germ using an affine conjugacy, setting the critical point at zero and the unique fixed point in $K(f_1)$ which does not cut $K(f_1)$ in two parts, the so-called $\beta$ fixed point of $f_1$, to 1). The operator $\mathcal{R}$ is called the Feigenbaum renormalization operator. In the setting of quadratic-like germs which have real values in the real line, there exists an unique fixed point to the Feigenbaum renormalization operator ([7]: see also [4]), denoted $f^\star$ (it is a open question if this is the unique fixed point in the set of all quadratic-like germs).

It is a consequence of the so-called a priori bounds [2] that we can choice a simply connected domain $U$, $\{0,1\} \subset U$, so that, if $\mathcal{B}(U)$ denotes the Banach space of the complex analytic functions $f$, $Df(0) = 0$, with a continuous extension to $\overline{U}$, provided with the sup norm, and $\mathcal{B}_{nor}(U)$ denotes the affine subspace of the functions $f$ so that $f(1) = 1$, then $f^\star$ has a complex analytic extension which belongs to $\mathcal{B}_{nor}(U)$ and, there exists $N$ so that the operator $\mathcal{R}^N$ can be represented

1991 *Mathematics Subject Classification.* 37F25, 37E20.
This work was supported by CNPq-Brazil grant 200764/01-2. I would like to thank M. Lyubich by my exciting visit to IMS-SUNY at Stony Brook.

as a compact operator defined in a small neighborhood of $f^\star$ in $\mathcal{B}_{nor}(U)$. More precisely, there exists a larger domain $\tilde{U} \Supset U$ and a complex analytic operator $\tilde{\mathcal{R}}\colon B_{\mathcal{B}_{nor}(U)}(f^\star,\epsilon) \to \mathcal{B}_{nor}(\tilde{U})$ so that, if $i$ denotes the natural inclusion $i\colon \mathcal{B}(\tilde{U}) \to \mathcal{B}(U)$, then $\mathcal{R}^N = i \circ \tilde{\mathcal{R}}$, where the equality holds in the intersection of the domains of the operators. To simplify the notation, we will assume that $N = 1$ and identify $\mathcal{R}$ with its complex analytic extension in $\mathcal{B}_{nor}(U)$.

Two quadratic-like maps $f$ and $g$ are in the same hybrid class if there exists a quasiconformal conjugacy $\phi$ between them, in a neighborhood of their filled-in Julia sets, so that $\overline{\partial}\phi \equiv 0$ on $K(f)$. Note that quadratic-like maps in the hybrid class of $f^\star$ are infinitly renormalizable. We will provide a new approach to the following result:

**Theorem 1.1** (Exponential contraction:[5] and [2]). *There exists $\lambda < 1$ so that, for every quadratic-like map $f$ which is in the same hybrid class that $f^\star$, there exist $n_0 = n_0(f)$ and $C = C(f) > 0$ so that $\mathcal{R}^n f \in B_{\mathcal{B}(U)}(f^\star,\epsilon)$, for $n \geq n_0$, and $|\mathcal{R}^{n_0+n} f - f^\star|_{\mathcal{B}(U)} \leq C\lambda^n$.*

The proof is more simple and straightforward than the previous methods, but indeed it is not necessary any breakthrough, so the reader should consider it as a natural development of (and interaction between) McMullen and Lyubich methods: see Remark 2 for more details. Other major attractive of this new proof is that it is essentially infinitesimal and has a "dynamical flavor": we will prove that the derivative of the renormalization operator is a contraction in the tangent space of the hybrid class. Moreover, the method seems to be so general as the previous ones: it also aplies to the classical renormalization horseshoe [2] and the Fibonacci renormalization operator [6], for instance.

## 2. Proof of the Theorem 1.1.

Let $f\colon V_1 \to V_2$ be a quadratic-like map with connected Julia set and with an analytical extension to $\mathcal{B}_{nor}(U)$, with $K(f) \subset U$. The horizontal subspace (introduced by Lyubich[2]) of $f$, denoted $E_f^h$, is the subspace of the vectors $v \in \mathcal{B}(U)$ so that there exists a quasiconfomal vector field in the Riemann sphere $\alpha$ satisfying $v = \alpha \circ f - Df \cdot \alpha$ in a neighborhood of $K(f)$, with $\overline{\partial}\alpha \equiv 0$ on $K(f)$ and $\alpha(0) = \alpha(1) = \alpha(\infty) = 0$. We will not use the following information here, but certainly it will clarify the spirit of our methods: in an apropriated setting, the hybrid class is a complex analytic manifold and the horizontal space is the tangent space of the hybrid class at $f$(see [2]).

**Lemma 2.1** ([1]). *Let $f$ be a quadratic-like map with an extension to $\mathcal{B}_{nor}(U)$ and connected Julia set contained in $U$. Assume that $f$ does not support invariant line fields in its filled-in Julia set. Let $V \Subset U$ be a domain with smooth boundary so that $K(f) \subset V$. Then there exist $C, \epsilon > 0$ so that, if $|f - g|_{\mathcal{B}(U)} \leq \epsilon$ and $g\colon g^{-1}V \to V$ is a quadratic-like map with connected Julia set, then, for every $v \in E_g^h$ there exists a $C|v|_{\mathcal{B}(U)}$-quasiconformal vector field $\alpha$ in $\overline{\mathbb{C}}$ so that $v = \alpha \circ g - Dg \cdot \alpha$ on $V$.*

With the aid of a compactness criterium to quasiconformal vector fields in $\overline{\mathbb{C}}$, we have:

**Corollary 2.2** ([1]). *Assume that $(f_n, v_n) \to_n (f_\infty, v_\infty)$ in $\mathcal{B}_{nor}(U) \times \mathcal{B}(U)$, where $f_i\colon f_i^{-1}V \to V$, $i \in \mathbb{N} \cup \{\infty\}$, are quadratic-like maps with connected filled-in Julia* sets $K(f_i) \subset V \Subset U$. Furthermore, assume that $v_n \in E^h_{f_n}$, for $n \in \mathbb{N}$. If $f$ does not support invariant line fields in $K(f)$, then $v_\infty \in E^h_f$. In particular $E^h_f$ is closed.

If $\mathcal{R}$ is the nth iteration of the Feigenbaum renormalization operator and $f$ is close to $f^\star$ in $\mathcal{B}(U)$, denote by $\beta_f$ the analytic continuation of the $\beta$-fixed point of the small Julia set associated with the nth renormalization of $f^\star$ [5]. The following result gives a description of the action of the derivative in a horizontal vector $v = \alpha \circ f - Df \cdot \alpha$ in terms of $\alpha$:

**Proposition 2.3.** *Let $V$ be a neighborhood of $K(f^\star)$. Replacing $\mathcal{R}$ by an iteration of it, if necessary, the following property holds: If $f \in \mathcal{B}_{nor}(U)$ is close enough to $f^\star$ and $v = \alpha \circ f - Df \cdot \alpha$ on $V$, where $v \in \mathcal{B}(U)$ and $\alpha$ is a quasiconformal vector field in the Riemann sphere, normalized by $\alpha(0) = \alpha(1) = \alpha(\infty) = 0$, then*

$$D\mathcal{R}_f \cdot v = r(\alpha) \circ \mathcal{R}f - D(\mathcal{R}f) \cdot r(\alpha), \tag{1}$$

*on $U$, where*

$$r(\alpha)(z) := \frac{1}{\beta_f}\alpha(\beta_f z) - \frac{1}{\beta_f}\alpha(\beta_f) \cdot z.$$

*In particular, if $f$ is renormalizable, then $D\mathcal{R}_f E^h_f \subset E^h_{\mathcal{R}f}$.*

This result is consequence of a simply calculation and the complex bounds [7] to $f^\star$. Note that, apart the normalization by a linear vector field, $r(\alpha)$ is just the pullback of the vector field $\alpha$ by a linear map. In particular, if $\alpha$ is a $C$-quasiconformal vector field, then $r(\alpha)$ is also a $C$-quasiconformal vector field: this will be a key point in the proof of the infinitesimal contraction of the renormalization operator in the horizontal subspace (Proposition 2.7).

Let $f^\star\colon V_1 \to V_2$ be a quadratic-like representation of the fixed point. The Feigenbaum tower is the indexed family of quadratic-like maps $f^\star_i\colon \beta_{f^\star}^{-i}V_1 \to \beta_{f^\star}^{-i}V_2$, $i \in \mathbb{N}$, defined by $f^\star_i(z) := \beta_{f^\star}^{-i}f^\star(\beta_{f^\star}^i \cdot z)$.

**Proposition 2.4** ([5])**.** *The Feigenbaum tower does not support invariant line fields: this means that there is not a measurable line field which is invariant by all (or even an infinite number of) maps in the Feigenbaum tower.*

**Proposition 2.5** ([7] and [5])**.** *Let $f$ be a quadratic-like map which admits a hybrid conjugacy $\phi$ with $f^\star$. Then $\phi_n(z) := \beta_{f^\star}^{-n} \cdot \phi(\beta_{\mathcal{R}^{n-1}f} \cdots \beta_f \cdot z)$ converges to identity uniformly on compact sets in the complex plane. In particular, there exists $n_0 = n_0(f)$ so that $\mathcal{R}^n f \in B_{\mathcal{B}_{nor}(U)}(f^\star, \epsilon)$, for $n > n_0$, and $\mathcal{R}^n f \to_n f^\star$ on $\mathcal{B}_{nor}(U)$.*

Theorem 1.1 says that this convergence is, indeed, exponentially fast. The following proposition has a straightforward proof:

**Proposition 2.6.** *Let $\mathcal{R} = i \circ \tilde{\mathcal{R}}$, where $\tilde{\mathcal{R}}\colon \mathcal{V} \subset \mathcal{B} \to \tilde{\mathcal{B}}$ is an operator defined in a neighborhood $\mathcal{V}$ of a Banach space $\mathcal{B}$, to another Banach space $\tilde{\mathcal{B}}$, and $i\colon \tilde{\mathcal{B}} \to \mathcal{B}$ is a compact linear transformation. Let $S \subset \mathcal{B} \times \mathcal{B}$ be a set with the following properties:*

1. **Vector bundle structure:** *If $(f, v_1)$ and $(f, v_2) \in S$, then $(f, \alpha \cdot v_1 + v_2) \in S$, for every $\alpha \in \mathbb{C}$,*
2. **Semicontinuity:** *If $(f_n, v_n) \to (f, v)$ and $(f_n, v_n) \in S$, then $(f, v) \in S$,*
3. **Invariance:** *If $(f, v) \in S$ then $(\mathcal{R}f, D\mathcal{R}_f \cdot v) \in S$,*
4. **Compactness:** *$\{(\tilde{\mathcal{R}}f, D\tilde{\mathcal{R}}_f \cdot v)\colon (f, v) \in S, |v| \leq 1\}$ is a bounded set in $\tilde{\mathcal{B}} \times \tilde{\mathcal{B}}$,*

5. **Uniform continuity:** Denote $E_f := \{(f, v) \colon (f, v) \in S\}$. There exists $C > 0$ so that $|D\mathcal{R}_f^n|_{E_f} \leq C$, for every $f$,
6. If $(f, v) \in S$ then $|D\mathcal{R}_f^n \cdot v| \to_n 0$,

Then there exist $\lambda < 1$ and $N \in \mathbb{N}$ so that $|D\mathcal{R}_f^N|_{E_f} \leq \lambda$, for every $f$ so that $E_f \neq \varnothing$.

Given $\epsilon$, $K > 0$, and a domain $V \Subset U$ so that $K(f^\star) \subset V$, denote by $\mathcal{A}(\epsilon, K, V)$ the set of maps $f \in \mathcal{B}_{nor}(U)$ so that there exists a $K$-quasiconformal map $\phi$ in the complex plane so that $\phi(V) \subset \overline{U}$ and $\phi \circ f^\star = f \circ \phi$ on $V$; moreover, for $n \geq 0$, we have $|\mathcal{R}^n f - f^\star|_{\mathcal{B}(U)} \leq \epsilon$. Note that $\mathcal{A} := \mathcal{A}(\epsilon, K, V)$ is closed. Furthermore, replacing $\mathcal{R}$ by an iteration, if necessary, we can assume that $\mathcal{A}$ is invariant by the action of $\mathcal{R}$. Selecting $K$ and $\epsilon$ properly, by the topological convergence (Proposition 2.5) and Lemma 2.2 in [3], for every $f$ in the hybrid class of $f^\star$, there exists $N = N(f)$ so that $\mathcal{R}^N f \in \mathcal{A}$.

**Proposition 2.7** (Infinitesimal contraction). *There exist $\lambda < 1$ and $N > 0$ so that $|D\mathcal{R}_f^N|_{E_f^h} \leq \lambda$, for every $f \in \mathcal{A}(\epsilon, K, V)$.*

*Proof.* Consider the set $S := \{(f, v) \colon f \in \mathcal{A}, \ v \in E_f^h\}$. It is sufficient to verify the properties in the statement of Proposition 2.6. Since $\mathcal{A}$ is closed, property 2 follows of Corollary 2.2. Since $\mathcal{A}$ is invariant by $\mathcal{R}$, property 3 follows of Propositon 2.3. The compactness property is obvious, if $\epsilon$ is small enough. To prove the uniform continuity property, by Propositions 2.1 and 2.3, we have that, for $(f, v) \in S$ and $n \geq 1$, $D\mathcal{R}_f^n \cdot v = \alpha_n \circ \mathcal{R}^N f - D(\mathcal{R}^N f) \cdot \alpha_n$ on $U$, with

$$\alpha_n(z) := \frac{1}{\beta_{n-1} \ldots \beta_0} \alpha(\beta_{n-1} \ldots \beta_0 z) - \frac{1}{\beta_{n-1} \ldots \beta_0} \alpha(\beta_{n-1} \ldots \beta_0) z,$$

where $\beta_i = \beta_{\mathcal{R}^i f}$ and $\alpha_n$ are $K \cdot |v|_{\mathcal{B}(U)}$-quasiconformal vector fields. Note that $K$ does not depends on $(f, v) \in S$ or $n \geq 1$. By the compactness of $K$-quasiconformal vector fields (recall that $\alpha_n(0) = \alpha_n(1) = \alpha_n(\infty) = 0$), we get $|D\mathcal{R}_f^n|_{E_f^h} \leq C$, for some $C > 0$. To prove assumption 6, note that $\overline{\partial} \alpha_n$ is an invariant Beltrami field to the finite tower

$$\mathcal{R}^n f, \frac{1}{\beta_{n-1}} \mathcal{R}^{n-1} f(\beta_{n-1} z), \ldots, \frac{1}{\beta_{n-1} \cdots \beta_0} f(\beta_{n-1} \cdots \beta_0 z).$$

But, by the topological convergence, these finite towers converges to the Feigenbaum tower. Hence, if a subsequence $\alpha_{n_k}$ converges to a quasiconformal vector field $\alpha_\infty$, then $\overline{\partial} \alpha_\infty$ is an invariant Beltrami field to the Feigenbaum tower (since $\overline{\partial} \alpha_{n_k}$ converges to $\overline{\partial} \alpha_\infty$ in the distributional sense), so, by Proposition 2.4, $\alpha_\infty$ is a conformal vector field in the Riemann sphere. Since $\alpha_\infty$ vanishes at three points, $\alpha_\infty \equiv 0$. Hence $\alpha_n \to 0$ uniformly on compact sets in the complex plane, so we get $\mathcal{R}_f^n \cdot v \to 0$ (Note that $|D(\mathcal{R}^n f)|$ is uniformly bounded, for $n \geq 1$). $\square$

*Remark* 1. *Proposition 2.6 is a generalization of the following fact about compact linear operators $T \colon \mathcal{B} \to \mathcal{B}$: if $T^n v \to 0$, for every $v \in \mathcal{B}$, then the spectral radius of $T$ is strictly smaller than one. Note that using only this simpler result, the non existence of invariant line fields in the Feigenbaum tower and Proposition 2.3 we can prove that the spectral radius of the operator $D\mathcal{R}_{f^\star}$ restricted to $E_{f^\star}^h$ is stricly smaller than one. Furthermore, to prove that the spectral radius is at most one, it is not even necessary the result about the Feigenbaum tower, but just the* *compactness of normalized quasiconformal vector fields and a well-known result about linear operators in Banach spaces: if $sup_n|T^n v| < \infty$ for every $v \in \mathcal{B}$, then $sup_n|T^n| < \infty$.*

We are going to prove Theorem 1.1: Let $f$ be a quadratic-like map in the hybrid class of $f^\star$. Then there exists a quasiconformal map $\phi\colon \mathbb{C} \to \mathbb{C}$ which is a conjugacy between them in a neighborhood of their Julia sets. Consider the following Beltrami path $f_t$ between the two maps, induced by $\phi$: if $\phi_t$, $|t| \leq 1$, is the unique normalized quasiconformal map so that $\overline{\partial}/\partial \phi_t = t \cdot \overline{\partial}/\partial \phi$, then $f_t = \phi_t \circ f \circ \phi_t^{-1}$. By the topological convergence, there exists $n_0$ so that $\mathcal{R}^{n_0+n} f_t \in \mathcal{A}$, for $n \geq 0$, $|t| \leq 1$. An easy calculation shows that

$$\left.\frac{d\mathcal{R}^{n_0+n} f_t}{dt}\right|_{t=t_0} \in E^h_{\mathcal{R}^{n_0+n} f_{t_0}},$$

for $|t_0| \leq 1$. The infinitesimal contraction finishes the proof.

*Remark 2.* There are two main steps in the above proof: the first one is to prove that $\alpha_n \to 0$ (in the proof of Proposition 2.7), which was inspired by a similar argument in a main technical result in McMullen proof: see the proof of Lemma 9.12 in [5]. The second step, to convert the topological convergence in exponential one, is almost straightforward here due to the compactness of the renormalization operator and a new result by Avila, Lyubich and de Melo, the semicontinuity of the horizontal directions [1]. In McMullen argument, additional considerations should be done to arrive in exponential contraction; firstly it is proved that quasiconformal deformations (as the quasiconformal vector field $\alpha$ in the definition of the horizontal vectors) are $C^{1+\beta}$-conformal at the critical point ( Lemma 9.12 in [5] and the deepness of the critical point have key roles in this proof), and then it is necessary to integrate this result. In Lyubich argument [2], firstly it is proved that the hybrid class is a complex analytic manifold and then the topological convergence is converted in exponential contraction via Schwartz's Lemma.

INSTITUTE FOR MATHEMATICAL SCIENCES, STATE UNIVERSITY OF NEW YORK AT STONY BROOK, STONY BROOK-NY, 11794-3660

*Home page*: **www.math.sunysb.edu/~smania/**
*E-mail address*: `smania@math.sunysb.edu`